\documentclass[draft]{birkjour}
\usepackage{latexsym,amssymb,mathrsfs,amsmath,bbm,eucal}%
\usepackage[mathx]{mathabx}
\allowdisplaybreaks

\def\noteq{\operatorname{=\hspace{-3.2mm}/\ }}
 \newtheorem{thm}{Theorem}[section]
 \newtheorem{cor}[thm]{Corollary}
 \newtheorem{lem}[thm]{Lemma}
 
 \theoremstyle{definition}
 \newtheorem{defn}[thm]{Definition}
 \theoremstyle{remark}

 \numberwithin{equation}{section}

\begin{document}

%-------------------------------------------------------------------------
% editorial commands: to be inserted by the editorial office
%
%\firstpage{1} \volume{228} \Copyrightyear{2004} \DOI{003-0001}
%\seriesextraline{This is the Concrete Title of this Book\br H.E. R and S.T.C. W, Eds.}
%\firstpage{1}

\title[The Hilbert-Schmidt analyticity]
 {The Hilbert-Schmidt analyticity associated with infinite-dimensional unitary groups}

%----------Author 1
\author[Lopushansky]{Oleh Lopushansky}

\address{%
1 Pigonia str.\\ 35-310 Rzesz\'ow\\
Poland}

\email{ovlopusz@ur.edu.pl}

\thanks{Faculty of Mathematics and Natural Sciences, Rzesz\'ow University.}



\date{January 11, 2015}

\begin{abstract}
The article is devoted to the problem of Hilbert-Schmidt type analytic extensions in Hardy spaces
over the infinite-dimensional unitary group endowed with an invariant probability measure.
Reproducing kernels of Hardy spaces, integral formulas of analytic extensions and their boundary values are considered.
\end{abstract}

%%% ----------------------------------------------------------------------
\maketitle
%%% ----------------------------------------------------------------------
%\tableofcontents
\section{Introduction}

The paper deals with the problem of Hilbert-Schmidt type analytic extensions
in the Hardy space  ${H}^2_\chi$ of complex  functions over the infinite-dimensional group
$U(\infty)=\bigcup\left\{ U(m)\colon m\in\mathbb{N}\right\}$  endowed with
an invariant probability measure $\chi$ where $U(m)$ are subgroups of unitary
 $m\times m$-matrices.
The measure $\chi$ is defined as a projective limit $\chi=\varprojlim\chi_m$ of the Haar probability
measures $\chi_m$ on $U(m)$. Moreover, $\chi$ is supported  by a projective limit $\mathfrak{U}=\varprojlim U(m)$ and
is invariant under the right action of $U^2(\infty):= {U(\infty)\times U(\infty)}$ on $\mathfrak{U}$.

A goal of this work is to find integral formulas for Hilbert-Schmidt analytic extensions
of functions from ${H}^2_\chi$ and to describe their radial  boundary values
on the open unit ball in a Hilbert space $\mathsf{E}$ where $U(\infty)$ acts irreducibly.

The measure $\chi$ on $\mathfrak{U}$ was described  by G.~\!Olshanski \cite{Olshanski2003},  Y.~\!Neretin \cite{Neretin2002}.
The notion $\mathfrak{U}$ is related to  D.~\!Pickrell's space of virtual Grassmannian \cite{Pickrell}.
Hardy spaces in infinite-dimensional settings were discussed in the works of
B.~\!Cole,  T.W.~\!Gamelin \cite{ColeGamelin86},  B.~\!{\O}rted, K.H.~\!Neeb \cite{OrtedNeeb98}.
Spaces of analytic functions of Hilbert-Schmidt holomorphy types were considered by T.A.W.~\!Dwyer III \cite{Dwyer71},
H.~\!Petersson \cite{Petersson2001}. More general classes of analytic functions associated with coherent sequences
of polynomial ideals were described by D.~\!Carando, V.~\!Dimant, S.~\!Muro \cite{Carado09}.
Integral formulas for analytic functions employing Wiener measures on infinite-dimensional Banach spaces
were suggested  by D.~\!Pinasco,  I.~\!Zalduendo  \cite{PinascoZalduendo05}.

Note that spaces of integrable functions with respect to invariant measures over
infinite-dimensional groups have been widely applied in stochastic processes
\cite{BorodinOlshanski05,Borodin11}, as well as in other areas.

This paper presents the following results.
In Theorem~\ref{irrep1}, we describe an orthogonal basis in the Hardy space ${H}^2_\chi$
indexed by means of Yang diagrams, consisting of $\chi$-essentially bounded functions.
Using this basis, in Theorem~\ref{Cauchy1} the reproducing kernel of ${H}^2_\chi$ is calculated.
It also allows us to define an antilinear isometric isomorphism $\mathcal{J}$ between ${H}^2_\chi$
and the symmetric Fock space $\Gamma$ generated by $\mathsf{E}$.
This isomorphism equips ${H}^2_\chi$ with a suitable infinite-dimensional analytic structure.
By means of $\mathcal{J}$, we establish in Theorem~\ref{hard3} an integral formula for
Hilbert-Schmidt analytic extensions of functions from ${H}^2_\chi$ on the open unit ball
$\mathsf{B}\subset\mathsf{E}$. The radial boundary values of these analytic extensions are described  in  Theorem~\ref{car:hardy2}.

\section{Background on invariant measure}

Let $U(m)$ $(m\in\mathbb{N})$ be the group of unitary $(m\times m)$-matrices. We endow $U(\infty)=\bigcup  U(m)$
with the inductive topology under every continuous inclusion $U(m)\looparrowright U(\infty)$
which assigns to any $u_m\in U(m)$ the matrix ${\begin{bmatrix}
                      u_m & 0\\
                      0 &\mathbbm{1}\\
                    \end{bmatrix}\in U(\infty)}$.
The right action over $U(\infty)$ is defined via
\begin{equation}\label{right0}
u.g=w^{-1}uv,\qquad u\in U(\infty),\quad g=(v, w)\in U^2(\infty)
\end{equation}
(the right action over $U(m)$ is defined similarly with ${u\in U(m)}$ and $g={(v, w)\in U^2(m)}$ where 
$U^2(m):= U(m)\times U(m)$).

Following \cite{Neretin2002,Olshanski2003},
every $u_m\in U(m)$ with $m>1$ can be written as
$u_m=\begin{bmatrix}
                      z_{m-1} & a \\
                      b & t \\
                    \end{bmatrix}$
so that $z_{m-1}$ is a $(m-1)\times(m-1)$-matrix and $t\in\mathbb{C}$.
It was proven that the Liv\v{s}ic-type mapping
(which is not a group homomorphism)
\begin{align}\label{projective}
&\pi^m_{m-1}\colon{u_m}\longmapsto u_{m-1}:=\left\{\!\!\begin{array}{lc}
                      z_{m-1}-[a(1+t)^{-1}b]&\!\!\!\!: t\not\eq-1 \\
                      z_{m-1}&\!\!\!\!: t=-1 \\
                    \end{array}\right.
\end{align}
from $U(m)$ onto $U(m-1)$ is Borel and surjective.

Consider the projective limit $\mathfrak{U}=\varprojlim U(m)$  taken with respect to $\pi_{m-1}^m$.
The embedding $\rho\colon U(\infty)\looparrowright\mathfrak{U}$
assigns to every  $u_m\in U(m)$  the stabilized sequence $u=(u_k)_{k\in\mathbb{N}}$
(see   \cite[n.4]{Olshanski2003}) so that
\begin{equation}\label{rho}
\rho\colon U(m)\ni u_m\longmapsto(u_k)\in\mathfrak{U},\qquad
u_k=\left\{\begin{array}{ccl}
\pi^m_k(u_m)&:&k<m,\\
                      u_m &:& k=m, \\
\begin{bmatrix}
                      u_m& 0 \\
                      0 &1\\
                    \end{bmatrix}
&:&k>m
\end{array}\right.
\end{equation}
where the projections $\pi_m\colon\mathfrak{U}\ni u\longrightarrow u_m\in U(m)$ such that ${\pi_{m-1}^m\circ\pi_m}=\pi_{m-1}$
are surjective and $\pi^m_k:=\pi^{k+1}_k\circ\ldots\circ \pi^m_{m-1}$ for ${k<m}$.
Using \eqref{right0}, the right action of $U^2(\infty)$ over  $\mathfrak{U}$ can be defined as
\begin{equation}\label{right}
\pi_m(u.g)=w^{-1}\pi_m(u)v,\qquad u\in \mathfrak{U}
\end{equation}
where $m$ is so large that $g=(v,w)\in U^2(m)$ (see  \cite[Def 4.5]{Olshanski2003}).

We endow every group $U(m)$ with the probability Haar measure $\chi_m$.
It is known \cite[Thm 1.6]{Neretin2002} that the pushforwards of $\chi_m$ to $U(m-1)$
under $\pi_{m-1}^m$ is the probability Haar measure $\chi_{m-1}$ on $U(m)$.
Let $U'(m)$ be the subset in $U(m)$ of matrices which do not have $\{-1\}$ as an eigenvalue.
Then $U'(m)$ is open in $U(m)$ and ${U(m)\setminus U'(m)}$ is $\chi_m$-negligible. Moreover,
the restriction $\pi_{m-1}^m\colon U'(m)\longrightarrow U'(m-1)$ is continuous and surjective
 \cite[Lem. 3.11]{Olshanski2003}.

Following \cite[Lem. 4.8]{Olshanski2003},  \cite[n.3.1]{Neretin2002}, via of the Kolmogorov consistency theorem
we uniquely define on $\mathfrak{U}$ the probability measure $\chi$ which is  the projective limit
under the mapping \eqref{projective}, i.e., we put
\begin{equation}\label{proj1}
\chi=\varprojlim\chi_m\quad\text{with}\quad \chi_m=\chi\circ\pi_m^{-1}
\quad\text{for all}\quad m\in\mathbb{N}.
\end{equation}
If $\mathfrak{U}'=\varprojlim U'(m)$ is the projective limit
with respect to $\pi_{m-1}^m\mid_{U'(m)}$ then
$\mathfrak{U}\setminus\mathfrak{U}'$ is $\chi$-negligible,
because $\chi_m$ is zero on $U(m)\setminus U'(m)$ for any $m$.

A complex-valued function on $\mathfrak{U}$ is called cylindrical if it has the form
$f=f_m\circ \pi_m$ for a certain $m\in\mathbb{N}$ and a complex function $f_m$
on $U(m)$ \cite[Def. 4.5]{Olshanski2003}.
By $L_\chi^\infty$ we denote the closed linear hull of all
cylindrical $\chi$-essentially bounded Borel  functions  endowed with  the norm
$\|f\|_{L_\chi^\infty}=\mathop{\rm ess\,sup}_{u\in\mathfrak{U}}|f(u)|$.

The measure \eqref{proj1} is a probability measure and is $U^2(\infty)$-invariant under the right actions
\eqref{right} over $\mathfrak{U}$ \cite[Prop. 3.2]{Neretin2002}.
Moreover, this measure is Radon so that
\begin{equation}\label{inv}
\int_\mathfrak{U}f(u.g)\,d\chi(u)=\int_\mathfrak{U}f(u)\,d\chi(u),\qquad g\in U^2(\infty),\quad f\in L_\chi^\infty
\end{equation}
and it satisfies  the property:  $(\chi\circ\pi_m^{-1})(K)=\chi_m(K)$ for any compact set $K$ in $U(m)$
 \cite[Lem. 1]{lopushansky2013}. Using the invariance property
\eqref{inv} and the Fubini theorem (see \cite[Lem. 2]{lopushansky2013}),  we obtain
\begin{align}\label{inv1}
\int_\mathfrak{U} f\,d\chi&=\int_\mathfrak{U}d\chi(u)\int_{U^2(m)}f(u.g)\,d(\chi_m\otimes\chi_m)(g),\\\label{inv2}
\int_\mathfrak{U}f\,d\chi&=\frac{1}{2\pi}\int_\mathfrak{U}\!d\chi(u)\int_{-\pi}^{\pi}f\left[\exp(\mathbbm{i}\vartheta)u\right]\,d\vartheta
\end{align}
for all $f\in L_\chi^\infty$. The closed linear hull of cylindrical complex functions  endowed with  the norm
$\|f\|_{L^2_\chi}=\left(\int_\mathfrak{U}|f|^2\,d\chi\right)^{1/2}$ is denoted by
$L^2_\chi$. It is clear that  $L^\infty_\chi\looparrowright L^2_\chi$ and
$\|f\|_{L^2_\chi}\le\|f\|_{L^\infty_\chi}$ for all ${f\in{L}^\infty_\chi}$.

\section{Hardy spaces}

Throughout the paper $\mathsf{E}$ is a separable complex Hilbert space with an orthonormal basis
$\left\{\mathfrak{e}_k\colon k\in\mathbb{N}\right\}$, scalar product $\langle\cdot\mid\cdot\rangle$ and
 norm $\|\cdot\|={\langle\cdot\mid\cdot\rangle^{1/2}}$.
So, for any element $x\in\mathsf{E}$  the following Fourier decomposition holds,
\begin{equation}\label{Fx}
x=\sum\mathfrak{e}_k\hat{x}_k,\qquad \hat{x}_k=\langle x\mid\mathfrak{e}_k\rangle.
\end{equation}
In what follows, let $\mathsf{B}=\left\{x\in\mathsf{E}\colon\|x\|<1\right\}$ and
$\mathsf{S}=\left\{x\in\mathsf{E}\colon\|x\|=1\right\}$.

Let $\mathsf{E}^{\otimes n}$ be the complete $n$th tensor power of
$\mathsf{E}$ endowed with the scalar product and norm
\[
\big\langle \psi\mid\phi\big\rangle=\langle x_1\mid y_1\rangle\ldots\langle x_n\mid y_n\rangle,
\qquad
\left\|\psi\right\|=\left\langle\psi\mid\psi\right\rangle^{1/2}
\]
for all  $\psi=x_1\otimes\ldots\otimes x_n$, $\phi=y_1\otimes\ldots\otimes y_n\in\mathsf{E}^{\otimes n}$
with $x_i,y_i\in\mathsf{E}$  ${(i=1,\ldots,n)}$.
As $\sigma\colon\{1,\ldots,n\}\longmapsto\{\sigma(1),\ldots,\sigma(n)\}$
runs through all $n$-elements permutations,
the symmetric complete $n$th tensor power $\mathsf{E}^{\odot n}$
is defined to be a codomain of the orthogonal projector
\[
\mathsf{E}^{\otimes n}\ni\psi\longmapsto
{x_1\odot\ldots\odot x_n}:= \frac{1}{n!}
\sum_\sigma{x_{\sigma(1)}\otimes\ldots\otimes x_{\sigma(n)}}\in\mathsf{E}^{\odot n}.
\]
Note that  $x^{\otimes n}={x\otimes\ldots\otimes x}={x\odot\ldots\odot x}=x^{\odot n}$.
Put $\mathsf{E}^{\otimes 0}=\mathsf{E}^{\odot 0}=\mathbb{C}$.

Let $\lambda=(\lambda_1,\ldots,\lambda_m)\in\mathbb{N}^m$
be a partition of an integer ${n\in\mathbb{N}}$ with ${m\le n}$
and ${\lambda_1\ge\lambda_2\ge\ldots\lambda_m>0}$, i.e.,
$|\lambda|=n$ where $|\lambda|:=\lambda_1+\ldots+\lambda_m$.
We identify partitions with Young diagrams.
By $\ell(\lambda)=m$ we denote the length of $\lambda$ defined
as  the number of rows in $\lambda$.  Let $\mathbb{Y}$ denote all Young diagrams and
$\mathbb{Y}_n:=\left\{\lambda\in\mathbb{Y}\colon|\lambda|=n\right\}$.
Assume that $\mathbb{Y}$ includes the empty partition  $\emptyset = (0, 0, \ldots )$.

An orthogonal basis in $\mathsf{E}^{\odot n}$ is formed by the system of symmetric tensor products
(see e.g. \cite[Sec. 2.2.2]{BerezanskiKondratiev95})
\[
\mathfrak{e}^{\odot\mathbb{Y}_n}=\bigcup_{\lambda\in\mathbb{Y}_n}
\left\{\mathfrak{e}^{\odot\lambda}_\imath:=
{\mathfrak{e}^{\otimes\lambda_1}_{\imath_1}\odot\ldots\odot\mathfrak{e}^{\otimes\lambda_m}_{\imath_m}}
\colon\imath\in\mathbb{N}^m_\ast, \ m=\ell(\lambda)\right\}, \quad \mathfrak{e}^{\odot\emptyset}_\imath=1
\]
where $\mathbb{N}^m_*:=
\left\{\imath=\left({\imath_1},\ldots,{\imath_m}\right)\in\mathbb{N}^m\colon\imath_j\noteq\imath_k,
\,\forall\,j\noteq k\right\}$. As is well known,
\begin{equation}\label{normfock}
\left\|\mathfrak{e}^{\odot\lambda}_\imath\right\|^2=
\frac{\lambda!}{n!},\qquad\lambda!:=\lambda_1!\cdot\ldots\cdot\lambda_m!.
\end{equation}

In what follows, we will use the fact that for every ${\psi\in\mathsf{E}^{\odot n}}$ one
can uniquely define the so-called {\it Hilbert-Schmidt $n$-homogenous polynomial}
\[\psi^*(x):={\left\langle x^{\otimes n}\mid\psi\right\rangle},\qquad {x\in\mathsf{E}}.\]
 In fact, the polarization formula  for symmetric tensor products (see \cite[1.5]{Floret97})
\begin{equation}\label{polarization}
z_1\odot\dots \odot z_n=\frac{1}{2^nn!}
\sum_{\theta_1,\ldots,\theta_n=\pm 1} \theta_1\dots \theta_n\,x^{\otimes n},\quad
x=\sum_{k=1}^n\theta_k z_k
\end{equation}
$(z_1,\ldots, z_n\in\mathsf{E})$ implies that the  $n$-homogenous polynomial
${\left\langle x^{\otimes n}\mid\psi\right\rangle}$ is uniquely defined by $\psi$, because
the set  ${z_1\odot\dots \odot z_n}$ is total in ${\mathsf{E}^{\odot n}}$.

Using the embedding \eqref{rho}, we define the $\mathsf{E}$-valued mapping
\[
\zeta\colon{\mathfrak{U}\ni u\longmapsto\rho^{-1}(u)\mathfrak{e}_1}
\]
which do not depend on the choice of $\mathfrak{e}_1$ in
\[\mathsf{S}(\infty):=\left\{\zeta(u)\colon u\in\mathfrak{U}\right\}=\bigcup\left\{\mathsf{S}(m)\colon m\in\mathbb{N}\right\}\]
where $\mathsf{S}(m)$ is the $m$-dimensional unit sphere.
In fact, for each stabilized sequence
$u=(u_k)\in\mathfrak{U}$ there exists an index $m$ such that $\rho^{-1}(u)\mathfrak{e}_1=u_k\mathfrak{e}_1$
belongs to  $\mathsf{S}(m)$  for all ${k\ge m}$.
On the other hand, for each $\mathfrak{e}\in\mathsf{S}(k)$ there exists ${v\in U(k)}$ such that
$v\mathfrak{e}=\mathfrak{e}_1$. Defining ${u.g\in\mathfrak{U}}$ with ${g=(1,v)\in U^2(k)}$
by means of \eqref{rho}-\eqref{right}, we have $\rho^{-1}(u.g)\mathfrak{e}=\pi_k(u.g)\mathfrak{e}=\pi_k(u)\mathfrak{e}_1
=\rho^{-1}(u)\mathfrak{e}_1$.

Consider the following system of cylindrical Borel functions
\[
\varepsilon_k(u):=
\left\langle\zeta(u)\mathrel{\big|}\mathfrak{e}_k\right\rangle,
\qquad k\in\mathbb{N}\]
where $\varepsilon_k:=\mathfrak{e}_k^*\circ\zeta$.
Using  $\zeta$, we may define the $\mathsf{E}^{\odot n}$-valued Borel  mapping
\[
\zeta^{\otimes n}\colon\mathfrak{U}\ni u\longmapsto
{\underbrace{\zeta(u)\otimes\ldots\otimes\zeta(u)}}_n,\qquad\zeta^{\otimes 0}\equiv1.\]
The following  assertion, which is a consequence of  the polarization formula \eqref{polarization},
is proved in  \cite[Lem. 3]{lopushansky2013}.

\begin{lem}\label{irrep}
The equality $\mathsf{S}(\infty)=\left\{\zeta(u)\colon u\in\mathfrak{U}'\right\}$ holds.
As a consequence,  to every  ${\psi\in\mathsf{E}^{\odot n}_\imath}$  there uniquely corresponds
the function in $L^\infty_\chi$
\[
\psi_\zeta(u):=\left\langle\zeta^{\otimes n}(u)\mathrel{\big|}\psi\right\rangle,\qquad
u\in\mathfrak{U}
\]
given by continuous restriction to $\mathfrak{U}'$. In particular,
to every ${\mathfrak{e}^{\odot\lambda}_\imath\in\mathfrak{e}^{\odot\mathbb{Y}_n}}$
there corresponds in $L^\infty_\chi$  the cylindrical function in the variable $u\in\mathfrak{U}$,
\begin{equation}\label{base2}
\varepsilon^{\lambda}_\imath(u)
:=\left\langle\zeta^{\otimes n}(u)\mathrel{\big|}
\mathfrak{e}^{\odot \lambda}_\imath\right\rangle=\prod_{k=1}^{\ell(\lambda)}
\left\langle\zeta(u)\mathrel{\big|}\mathfrak{e}_{\imath_k}\right\rangle^{\lambda_k}.
\end{equation}
\end{lem}

Lemma~\ref{irrep} straightforwardly implies that the system
$\mathfrak{e}^{\odot\mathbb{Y}}:=\bigcup\mathfrak{e}^{\odot\mathbb{Y}_n}$
of tensor products $\mathfrak{e}^{\odot\lambda}_\imath=
\mathfrak{e}^{\otimes\lambda_1}_{\imath_1}\odot\ldots\odot\mathfrak{e}^{\otimes\lambda_m}_{\imath_m}$,
indexed by $\lambda={(\lambda_1,\ldots,\lambda_m)\in\mathbb{Y}}$ and
$\imath={\left({\imath_1},\ldots,{\imath_m}\right)\in\mathbb{N}^m_\ast}$ with $m=\ell(\lambda)$,
uniquely defines the appropriate system
\[
\varepsilon^\mathbb{Y}:=\bigcup_{\lambda\in\mathbb{Y}}
\left\{\varepsilon^{\lambda}_\imath:=
\varepsilon^{\lambda_1}_{\imath_1}\odot\ldots\odot\varepsilon^{\lambda_m}_{\imath_m}
\colon\imath\in\mathbb{N}^m_\ast,  \ m=\ell(\lambda)\right\},\quad \varepsilon^\emptyset_\imath\equiv1,
\]
of $\chi$-essentially bounded cylindrical functions in the variable  $u\in\mathfrak{U}$ that
possess continuous restrictions to $\mathfrak{U}'$.

\begin{thm}\label{irrep1}
For any $\imath\in\mathbb{N}^m_\ast$ and $\psi,\phi\in\mathsf{E}^{\odot n}_\imath$,  the following equality holds,
\begin{equation}\label{polin}
\binom{n+m-1}{n}\int_\mathfrak{U}\phi_\zeta\,\bar\psi_\zeta\,d\chi=\left\langle\psi\mid \phi\right\rangle.
\end{equation}
As a consequence, given $(\lambda,\imath)\in\mathbb{Y}\times\mathbb{N}^m_\ast$ with $m=\ell(\lambda)$,
the system $\varepsilon^\mathbb{Y}$ of functions $\varepsilon^\lambda_\imath$ is orthogonal in the space $L_\chi^2$ and
\begin{equation}\label{norm}
\left\|\varepsilon^{\lambda}_\imath\right\|_{L^2_\chi}=\left(\frac{(m-1)! \lambda!}{(m-1+|\lambda|)!}\right)^{1/2}.
\end{equation}
\end{thm}

\begin{proof}
Let $\mathsf{E}_\imath$ with $\imath=\left(\imath_1,\ldots,\imath_m\right)\in\mathbb{N}^m_\ast$
be the $m$-dimensional subspace in $\mathsf{E}$ spanned by
$\left\{\mathfrak{e}_{\imath_1},\ldots,\mathfrak{e}_{\imath_m}\right\}$ and $U(\imath)$ be the unitary subgroup of
$U(\infty)$ acting in $\mathsf{E}_\imath$. The symbol
$\mathsf{E}^{\odot n}_\imath$ means the $n$th symmetric tensor power
of $\mathsf{E}_\imath$. Briefly denote $\psi_\dag[v\zeta(u)]:=
\big\langle\big([v\rho^{-1}(u)]\mathfrak{e}_1\big)^{\otimes n}\mathop{\big|}\psi\big\rangle$
with ${\psi\in\mathsf{E}^{\odot n}_\imath}$
for all ${v\in U(\imath)}$ and ${u\in\mathfrak{U}}$. Using \eqref{inv1}
with $U(\imath)$ instead of $U(m)$, we have
\begin{equation}\label{iso0}
\int_\mathfrak{U}\phi_\zeta\,\bar\psi_\zeta\,d\chi=\int_\mathfrak{U}d\chi(u)
\int_{U(\imath)}\phi_\dag[v\zeta(u)]\cdot\bar\psi_\dag[v\zeta(u)]\,d\chi_\imath(v)
\end{equation}
for all $\psi,\phi\in\mathsf{E}^{\odot n}_\imath$. It is clear that
\[
\Big|\int_{U(\imath)}\phi_\dag\,\bar\psi_\dag\,d\chi_\imath\Big|\le
\sup_{v\in {U(\imath)}}\big|\phi_\dag[v\zeta(u)]\big|\,
\big|\psi_\dag[v\zeta(u)]\rangle\big|\le\|\phi\|\,\|\psi\|
\]
for all ${u\in\mathfrak{U}}$.
Hence, the corresponding sesquilinear form in \eqref{iso0} is continuous on $\mathsf{E}^{\odot n}_\imath$.
Thus, there exists a linear bounded operator $A$ over $\mathsf{E}^{\odot n}_\imath$ such that
\[\left\langle A\psi\mid\phi\right\rangle=\int_{U(\imath)}\phi_\dag\,\bar\psi_\dag\,d\chi_\imath.\]

Next we show that  $A$ commutes with all operators
$w^{\otimes n}\in\mathscr{L}\left(\mathsf{E}^{\odot n}_\imath\right)$ with $w\in {U(\imath)}$
acting as $w^{\otimes n}x^{\otimes n}=(wx)^{\otimes n}$, ${(x\in\mathsf{E}_\imath)}$.
Invariant properties \eqref{inv} of $\chi_\imath$ under the right action \eqref{right} yield
\[\begin{split}
&\left\langle(A\circ w^{\otimes n})\psi\mid\phi\right\rangle=\\
&=\int_{U(\imath)}\left\langle[v\zeta(u)]^{\otimes n}\mid\phi\right\rangle
\overline{\left\langle[v\zeta(u)]^{\otimes n}\mid w^{\otimes n}\psi\right\rangle}d\chi_\imath(v)\\
&=\int_{U(\imath)}\left\langle[w^{-1}v\zeta(u)]^{\otimes n}
\mid (w^{-1})^{\otimes n}\phi\right\rangle
\overline{\left\langle[w^{-1}v\zeta(u)]^{\otimes n}
\mid \psi\right\rangle}d\chi_\imath(v)\\
&=\int_{U(\imath)}\left\langle[v\zeta(u)]^{\otimes n}\mid (w^{-1})^{\otimes n}\phi\right\rangle
\overline{\left\langle[v\zeta(u)]^{\otimes n}\mid \psi\right\rangle}d\chi_\imath(v)\\
&=\left\langle A\psi\mid(w^{-1})^{\otimes n}\phi\right\rangle
=\left\langle(w^{\otimes n}\circ A)\psi\mid\phi\right\rangle,
\end{split}\]
where  $w^{-1}\in {U(\imath)}$ is the hermitian adjoint matrix of $w$. Hence, the equality
\begin{equation}\label{schur}
A\circ w^{\otimes n}=w^{\otimes n}\circ A,\qquad {w\in U(\imath)}
\end{equation}
holds. Let us check that the operator $A$, satisfying the condition \eqref{schur},
 is proportional to the identity operator on $\mathsf{E}^{\otimes n}_\imath$.
To this end we form the $n$th tensor power of the unitary group $U(\imath)$,
\[[U(\imath)]^{\otimes n}=\left\{
w^{\otimes n}\in\mathscr{L}\left(\mathsf{E}^{\odot n}_\imath\right)\colon w\in U(\imath)\right\},
\qquad [U(\imath)]^{\otimes 0}=1.\]
Clearly, $[U(\imath)]^{\otimes n}$ is a unitary group over $\mathsf{E}^{\odot n}_\imath$. Let us
check that the corresponding unitary representation
\begin{equation}\label{diag}
U(\imath)\ni w\longmapsto w^{\otimes n}\in\mathscr{L}\left(\mathsf{E}^{\odot n}_\imath\right)
\end{equation}
is irreducible. This means that  there is no subspace in $\mathsf{E}^{\odot n}_\imath$
other than $\{0\}$ and the whole space which is invariant under the action of $[U(\imath)]^{\otimes n}$.

Suppose, on the contrary, that there is an element ${\psi\in\mathsf{E}^{\odot n}_\imath}$ such that the equality
${\big\langle\big([w\rho^{-1}(u)]\mathfrak{e}_1\big)^{\otimes n}\mathop{\big|}\psi\big\rangle=0}$
holds for all ${w\in U(\imath)}$ and ${u\in U(\infty)}$. By Lemma~\ref{irrep} the
elements $w\rho^{-1}(u)$ act transitively on  $\mathsf{S}(\infty)$. Hence, by $n$-homogeneity, we obtain
${\langle x^{\otimes n}\mid\psi\rangle=0}$ for all ${x\in\mathsf{E}_\imath}$. Applying the polarization formula \eqref{polarization},
we get ${\psi=0}$. Hence, \eqref{diag} is irreducible.

Thus, we can apply to \eqref{diag}  the Schur lemma
\cite[Thm 21.30]{HewittRoss70}: a non-zero matrix which commutes with all matrices of an irreducible representation is a
constant multiple of the unit matrix. As a result,
we obtain that the operator $A$, satisfying  \eqref{schur},  is proportional to
the identity operator on $\mathsf{E}^{\odot n}_\imath$ i.e.
$A =\alpha_{(n,\imath)}\mathbbm{1}_{\mathsf{E}^{\odot n}_\imath}$ with a constant  $\alpha_{(n,\imath)}>0$.
It follows that
\begin{equation}\label{MainEq}
\int_{U(\imath)}\phi_\dag\,\bar\psi_\dag\,d\chi_\imath=
\alpha_{(n,\imath)}\left\langle \psi\mid\phi\right\rangle,\qquad
\phi,\psi\in\mathsf{E}^{\odot n}_\imath.
\end{equation}
In particular, the subsystem of cylindrical functions $\varepsilon^\lambda_\imath$ with a fixed ${\imath\in\mathbb{N}^m_\ast}$
is orthogonal in $L_\chi^2$, because the corresponding system of tensor products $\mathfrak{e}^{\odot \lambda}_\imath$
indexed by $\lambda\in\mathbb{Y}_n$ with ${\ell(\lambda)=m}$ forms an orthogonal basis in $\mathsf{E}^{\odot n}_\imath$.

It remains to note that  the set of all indices $\imath={\left({\imath_1},\ldots,{\imath_m}\right)\in\mathbb{N}^m_\ast}$
with all $m=\ell(\lambda)$  is directed with respect to  the set-theoretic embedding,  i.e., for any $\imath,\imath'$
there exists $\imath''$ so that $\imath\cup\imath'\subset\imath''$. This fact and the above reasoning imply
that the whole system $\varepsilon^\mathbb{Y}$ is also orthogonal in $L_\chi^2$.

Taking into account \eqref{normfock},
we can choose $\phi_n=\psi_n=\varepsilon^{\lambda}_\imath\sqrt{n!/\lambda!}$  in \eqref{MainEq}.
As a result, we obtain
\[
\alpha_{(n,\imath)}=\frac{n!}{\lambda !}\int_{U(\imath)}\left|\varepsilon^\lambda_\imath\right|^2d\chi_\imath
=\frac{n!}{\lambda !}\left\|\varepsilon^\lambda_\imath\right\|^2_{L^2_\chi}.
\]
The well known formula \cite[1.4.9]{RudinFT80} for the unitary $m$-dimensional group gives
\[
\int_{U(\imath)}\left|\varepsilon^\lambda_\imath\right|^2d\chi_\imath=
\frac{\lambda !(m-1)!}{(n+m-1)!},\qquad|\lambda|=n,\quad{\ell(\lambda)=m}.
\]
Using the last two formulas, we arrive at the relation
\begin{equation}\label{rudin}
\alpha_{(n,\imath)}=
\frac{n!}{\lambda !}\int_{U(\imath)}\left|\varepsilon^\lambda_\imath\right|^2d\chi_\imath
=\frac{n!}{\lambda !}\,\frac{\lambda !(m-1)!}{(n+m-1)!}=\frac{n!(m-1)!}{(n+m-1)!}.
\end{equation}
Combining \eqref{iso0} and \eqref{rudin}, we get  \eqref{polin} and, as a consequence, \eqref{norm}.
\end{proof}

\begin{defn}
By $H_\chi^2$ we  denote  the Hardy space over  $U(\infty)$ defined as
the $L^2_\chi$-closure of the complex linear span of the orthogonal system $\varepsilon^\mathbb{Y}$.
\end{defn}

Let the space $H_\chi^{2,n}$  be the $L^2_\chi$-closure of the complex linear span of the subsystem
$\varepsilon^{\mathbb{Y}_n}:=
\big\{\varepsilon^\lambda_\imath\in\varepsilon^\mathbb{Y}\colon{(\lambda,\imath)\in\mathbb{Y}_n\times\mathbb{N}^{\ell(\lambda)}_\ast}\big\}$
with a fixed ${n\in\mathbb{Z}_+}$.

\begin{cor}\label{ortog}
For any positive integers $n\noteq k$  the orthogonality
$H_\chi^{2,n}\perp H_\chi^{2,k}$ in $L^2_\chi$ holds.  As a consequence, the following orthogonal decomposition holds,
\begin{equation}\label{ort}
H_\chi^2=\mathbb{C}\oplus H_\chi^{2,1}\oplus H_\chi^{2,2}\oplus\ldots.
\end{equation}
\end{cor}
\begin{proof}
The orthogonal property ${\varepsilon^\mu_\jmath\perp\varepsilon^{\lambda}_\imath}$ with
${|\mu|\not\eq|\lambda|}$ for any ${\imath\in\mathbb{N}^{\ell(\lambda)}_\ast}$ and
${\jmath\in\mathbb{N}^{\ell(\mu)}_\ast}$  follows from \eqref{inv2}, since
\[\begin{split}
\int_\mathfrak{U}\varepsilon^\mu_\jmath\,\bar\varepsilon^\lambda_\imath\,d\chi&=
\int_\mathfrak{U}
\varepsilon^\mu_\jmath\big(\exp(\mathbbm{i}\vartheta)u\big)\,
\bar\varepsilon^\lambda_\imath\big(\exp(\mathbbm{i}\vartheta)u\big)d\chi(u)\\
&=\frac{1}{2\pi}\int_\mathfrak{U}\varepsilon^\mu_\jmath\,\bar\varepsilon^\lambda_\imath\,d\chi
\int_{-\pi}^\pi{\exp\big(\mathbbm{i}(|\mu|-|\lambda|)\vartheta\big)}\,d\vartheta=0
\end{split}\]
for all  $\lambda\in\mathbb{Y}$ and $\mu\in\mathbb{Y}\setminus\{\emptyset\}$.
This yields  $H_\chi^{2,|\mu|}\perp H_\chi^{2,|\lambda|}$ in the space $L^2_\chi$.
\end{proof}

\section{Reproducing kernels}

Let us construct the reproducing kernel of $H_\chi^2$. We refer to \cite{Sait} regarding reproducing kernels.

\begin{lem}\label{ReprodP}
For every $u,v\in\mathfrak{U}$ there exists ${q\in\mathbb{N}}$ such that the reproducing kernel of the subspace $H_\chi^{2,n}$ in $L^2_\chi$
has the form
\begin{equation}\label{reprod}
\begin{split}
\mathfrak{h}_n(v,u)&=\sum_{m\le q}
\binom{n+m-1}{n}
\left\langle\zeta(v)\mid\zeta(u)\right\rangle^n\\
&=\sum_{(\lambda,\imath)\in\mathbb{Y}_n\times\mathbb{N}^{\ell(\lambda)}_\ast}
\frac{\varepsilon^\lambda_\imath(v)\,\bar\varepsilon^\lambda_\imath(u)}
{\|\varepsilon^{\lambda}_\imath\|^2_{L^2_\chi}}, \qquad u,v\in\mathfrak{U}.
\end{split}
\end{equation}
\end{lem}
\begin{proof}
Note that $\mathfrak{h}_0\equiv1$. From \eqref{rho} it follows that for each stabilized sequence
${u\in\mathfrak{U}}$ there exists ${u_m\in U(m)}$ with a certain $m=m(u)$ such that $u=\rho(u_m)$.
So, the element $\zeta(u)=\rho^{-1}(u)\mathfrak{e}_1$
is located on the $m$-dimensional sphere $\mathsf{S}(m)$.
It means that its Fourier series $\zeta(u)={\sum{\mathfrak e}_k\varepsilon_k(u)}$ has $m(u)$ terms.
The tensor multinomial theorem yields the  Fourier decomposition
\[
[\zeta(u)]^{\otimes n}=
\left(\sum{\mathfrak e}_k\varepsilon_k(u)\right)^{\otimes n}=\sum_{(\lambda,\imath)\in\mathbb{Y}_n\times\mathbb{N}^{\ell(\lambda)}_\ast}
\frac{n!}{\lambda!}\mathfrak{e}^{\odot\lambda}_\imath\,\varepsilon^{\lambda}_\imath(u)
\]
in the space $\mathsf{E}^{\odot n}$. Using the formula \eqref{normfock}, we obtain
\begin{align*}
&\left\langle \zeta(v)\mid\zeta(u)\right\rangle^n=
\left\langle [\zeta(v)]^{\otimes n}\mid[\zeta(u)]^{\otimes n}\right\rangle\\
&=\sum_{(\lambda,\imath)\in\mathbb{Y}_n\times\mathbb{N}^{\ell(\lambda)}_\ast}\Big(\frac{n!}{\lambda!}\Big)^2
\left\langle \mathfrak{e}^{\odot\lambda}_\imath\mid\mathfrak{e}^{\odot\lambda}_\imath\right\rangle
\varepsilon^\lambda_\imath(v)\,\bar\varepsilon^\lambda_\imath(u)
=\sum_{(\lambda,\imath)\in\mathbb{Y}_n\times\mathbb{N}^{\ell(\lambda)}_\ast}
\frac{\varepsilon^\lambda_\imath(v)\,\bar\varepsilon^\lambda_\imath(u)}
{\|\varepsilon^{\lambda}_\imath\|^2_{L^2_\chi}}
\end{align*}
where  $\left\langle \zeta(v)\mid\zeta(u)\right\rangle$ is decomposed into $q=\min\{m(u),m(v)\}$ summands in virtue of orthogonality.
Multiplying both sides by $\binom{n+m-1}{n}$ and summing over all $m\le q$,  we get  \eqref{reprod}. It follows that
$\int_\mathfrak{U}\mathfrak{h}_n(v,u)\varepsilon^\lambda_\imath(u)\,d\chi(u)=\varepsilon^\lambda_\imath(v)$ for each ${v\in\mathfrak{U}}$.
Via Theorem~\ref{irrep} the system $\varepsilon^{\mathbb{Y}_n}$ of functions $\varepsilon^\lambda_\imath$
 forms an orthogonal basis in  $H_\chi^{2,n}$.  So, the integral operator
\begin{equation}\label{Taypol}
\int_\mathfrak{U}\mathfrak{h}_n(v,u)\psi_\zeta(u)\,d\chi(u)=\psi_\zeta(v),\qquad {\psi_\zeta\in H_\chi^{2,n}}
\end{equation}
acts identically on $H_\chi^{2,n}$.
Thus, the kernel \eqref{reprod} is reproducing in $H_\chi^{2,n}$.
\end{proof}

Let us consider the complex-valued kernel
\[
\mathfrak{h}(z;v,u)=\prod_{m\le\min\{m(u),m(v)\}}
\left[{\phantom{\big|}}\!\!1-z\left\langle\zeta(v)\mid\zeta(u)\right\rangle\right]^{-m},\quad
u,v\in\mathfrak{U},\quad |z|<1
\]
where $m(u)$ is the number of terms in the Fourier series $\zeta(u)={\sum{\mathfrak e}_k\varepsilon_k(u)}$.

\begin{thm}\label{Cauchy1}
The expansion $\mathfrak{h}(z;v,u)=\sum z^n\mathfrak{h}_n(v,u)$ holds for any ${u,v\in\mathfrak{U}}$ and ${|z|<1}$.
The kernel $\mathfrak{h}(1;v,u)=\sum\mathfrak{h}_n(v,u)$ is reproducing in $H^2_\chi$ in the sense that
\begin{equation}\label{sum}
\int_\mathfrak{U}\mathfrak{h}(1;v,u)f(u)\,d\chi(u)=f(v),\qquad{f\in H_\chi^2},\quad v\in\mathfrak{U}.
\end{equation}
\end{thm}

\begin{proof}
Let $q=\min\{m(u),m(v)\}$ and $m\le q$. As is well known \cite[1.4.10]{RudinFT80},
\begin{equation}\label{prod}
\left[{\phantom{\big|}}\!\!1- z\left\langle\zeta(v)\mid
\zeta(u)\right\rangle\right]^{-m}=\sum_{n\in\mathbb{Z}_+}\binom{n+m-1}{n}
\left\langle z\zeta(v)\mid\zeta(u)\right\rangle^n
\end{equation}
for all ${|z|<1}$. By the Vandermonde identity, we have
\begin{align*}
&\binom{n+m-1}{n}\left\langle z\zeta(v)\mid\zeta(u)\right\rangle^n=
\binom{r+k+p+l-2}{r+k}\left\langle z\zeta(v)\mid\zeta(u)\right\rangle^{r+k}\\
&=\sum_{r=0}^n\binom{r+p-1}{r}\binom{n-r+l-1}{n-r}\left\langle z\zeta(v)\mid\zeta(u)\right\rangle^{r+k}
\end{align*}
for all $n=r+k$ and $m=p+l-1$. Applying recursively this identity to  the series \eqref{prod} with any ${m\le q}$
and using Lemma~\ref{ReprodP}, we obtain
\begin{align*}
\mathfrak{h}(z;v,u)&=\prod_{m\le q}\sum_{n\in\mathbb{Z}_+}
\binom{n+m-1}{n}\left\langle z\zeta(v)\mid\zeta(u)\right\rangle^n\\
&=\sum_{n\in\mathbb{Z}_+}z^n
\sum_{(\lambda,\imath)\in\mathbb{Y}_n\times\mathbb{N}^{\ell(\lambda)}_\ast}
\frac{\varepsilon^\lambda_\imath(v)\,\bar\varepsilon^\lambda_\imath(u)}
{\|\varepsilon^{\lambda}_\imath\|^2_{L^2_\chi}}=\sum_{n\in\mathbb{Z_+}}z^n\mathfrak{h}_n(v,u).
\end{align*}
Hence, the required expansion holds. By  \eqref{ort} we have $f=\sum_nf_n$ for any ${f\in H_\chi^2}$ where
$f_n\in H_\chi^{2,n}$ is the orthogonal projection of $f$.
Observing that  $\mathfrak{h}_k(z;\cdot,u)\perp f_n(\cdot)$ with $n\noteq k$ holds in $L^2_\chi$, we obtain
\[
\int_\mathfrak{U}\mathfrak{h}(1;v,u)f(u)\,d\chi(u)=
\sum \int_\mathfrak{U}\mathfrak{h}_n(v,u)f_n(v)\,d\chi(u)=\sum  f_n(v)=f(v)
\]
for all $v\in\mathfrak{U}$ and ${f\in H_\chi^2}$. Hence, \eqref{sum} is valid.
\end{proof}

\section{The Hilbert-Schmidt analyticity}

Recall (see e.g. \cite{G}) that a function $f$ on an open domain  in a Banach space
is said to be analytic if it is G\^ateaux analytic and norm continuous.
Similarly to \cite{Dwyer71,Petersson2001}, we say that $f$ is {\it Hilbert-Schmidt analytic}
if its Taylor coefficients are Hilbert-Schmidt polynomials.
Now we describe a space $H^2$ of Hilbert-Schmidt analytic complex functions on the open
ball $\mathsf{B}$.

The symmetric Fock space is defined to be the orthogonal sum
\[
\Gamma=\bigoplus_{n\in\mathbb{Z}_+}\mathsf{E}^{\odot n},\qquad
\langle \psi\mid\phi\rangle=\sum_{n\in\mathbb{Z}_+}\langle \psi_n\mid\phi_n\rangle
\]
for all elements $\psi=\bigoplus_n\psi_n$, $\phi=\bigoplus_n\phi_n\in\Gamma$ with ${\psi_n,\phi_n\in\mathsf{E}^{\odot n}}$.
The subset $\left\{x^{\otimes n}\colon x\in\mathsf{B}\right\}$ is total in $\mathsf{E}^{\odot n}$
by virtue of \eqref{polarization}. This provides the total property of the subsets
$\left\{(1-x)^{-\otimes1}\colon x\in\mathsf{B}\right\}$ in  $\Gamma$ where we denote
\[
(1-x)^{-\otimes1}:=\sum x^{\otimes n},\qquad x^{\otimes 0}=1.
\]
The $\Gamma$-valued function $(1-x)^{-\otimes1}$ in the variable $x\in\mathsf{B}$ is analytic, since
\begin{equation}\label{ob}
\left\|(1-x)^{-\otimes1}\right\|^2=\sum\|x\|^{2n}=\left(1-\|x\|^2\right)^{-1}<\infty.
\end{equation}

Let us define the Hilbert space of analytic complex functions
 in the variable $x\in\mathsf{B}$, associated with the Fock space $\Gamma$, as follows
$$
H^2=\left\{\psi^*(x)=\left\langle(1-x)^{-\otimes1}\mid\psi\right\rangle\colon \psi\in\Gamma\right\},\qquad
\left\|\psi^*\right\|_{H^2}:=\left\|\psi\right\|
$$
for all  $x\in\mathsf{B}$. This description is correct, because each  function $\psi^*$ in the variable $x\in\mathsf{B}$
is analytic by virtue of \cite[Prop. 2.4.2]{Herv},
as a composition of the  analytic $\Gamma$-valued function $(1-x)^{-\otimes1}$
in the variable ${x\in\mathsf{B}}$ and the linear functional $\left\langle\cdot\mid\psi\right\rangle$ on $\Gamma$.

Similarly, we define the closed subspace in $H^2$ of $n$-homogenous Hilbert-Schmidt polynomials $\psi_n^*$
in the variable $x\in\mathsf{E}$ as
$$
H^2_n=\left\{\psi_n^*(x)=\left\langle x^{\otimes n}\mid\psi_n\right\rangle\colon \psi_n\in\mathsf{E}^{\odot n}\right\}.
$$
Differentiating at zero any function $\psi^*={\bigoplus\psi^*_n\in H^2}$ with ${\psi^*_n\in H^2_n}$,
we obtain that its Taylor coefficients at zero $(n!)^{-1}d^n_0\psi^*=\psi^*_n$ are Hilbert-Schmidt polynomials.
Hence, every function from $H^2$ is Hilbert-Schmidt analytic. Clearly, the following orthogonal decomposition holds,
\begin{equation}\label{iso}
H^2=\mathbb{C}\oplus H_1^2\oplus H_2^2\oplus\ldots.
\end{equation}

One can show that $\left(H^2_n\right)_n$ is a coherent sequence of polynomial ideals over $\mathsf{E}$
in the meaning of \cite[Def. 1.1]{Carado09}.

For each pair ${(\lambda,\imath)\in\mathbb{Y}_n\times\mathbb{N}^{\ell(\lambda)}_\ast}$,
we can uniquely assign the Hilbert-Schmidt $n$-homogenous polynomial
\[
\hat{x}^\lambda_\imath:=\left\langle x^{\otimes n}\mathrel{\big|}\mathfrak{e}^{\odot \lambda}_\imath\right\rangle,\qquad x\in\mathsf{E},
\]
defined via the Fourier coefficients $\hat{x}_k:=\mathfrak{e}_k^*(x)={\langle x\mid\mathfrak{e}_k\rangle}$
of an element ${x\in\mathsf{E}}$. Taking into account \eqref{normfock}, the tensor multinomial theorem yields the following
orthogonal decompositions with respect to the  basis  $\mathfrak{e}^{\odot\mathbb{Y}}$ in $\Gamma$,
\begin{equation}\label{Tayl}
(1-x)^{-\otimes1}=
\sum_{(\lambda,\imath)\in\mathbb{Y}\times\mathbb{N}^{\ell(\lambda)}_\ast}
\frac{\hat{x}^\lambda_\imath\mathfrak{e}^{\odot\lambda}_\imath}{\|\mathfrak{e}^{\odot\lambda}_\imath\|^2},\qquad x\in\mathsf{B}.
\end{equation}
Hence, any function ${\psi^*\in H^2}$ has the orthogonal expansion
\begin{equation}\label{Pn1}
\psi^*(x)=
\left\langle(1-x)^{-\otimes1}\mid\psi\right\rangle=\sum_{(\lambda,\imath)\in\mathbb{Y}\times\mathbb{N}^{\ell(\lambda)}_\ast}
\hat\psi_{(\lambda,\imath)}{\hat{x}^\lambda_\imath},\qquad x\in\mathsf{B}
\end{equation}
where $\hat\psi_{(\lambda,\imath)}:=\langle\mathfrak{e}^{\odot\lambda}_\imath\mid\psi\rangle
\|\mathfrak{e}^{\odot\lambda}_\imath\|^{-2}$ are the Fourier coefficients of ${\psi\in\Gamma}$
with respect to the basis $\mathfrak{e}^{\odot\mathbb{Y}}$  and, moreover,
$\|\psi^*\|_{H^2}^2=\sum_{(\lambda,\imath)}
|\langle\mathfrak{e}^{\odot\lambda}_\imath\mid\psi\rangle|^2\|\mathfrak{e}^{\odot\lambda}_\imath\|^{-2}$.
Thus, $\|\psi^*\|_{H^2}$ is a Hilbert-Schmidt type norm on $H^2$.

\section{Integral formulas}

The one-to-one correspondence $\mathfrak{e}^{\odot \lambda}_\imath\leftrightarrow\varepsilon^\lambda_\imath$ allows us to construct
an antilinear isometric isomorphism $\mathcal{J}\colon\Gamma\longrightarrow{H}^2_\chi$
and its adjoint   ${\mathcal{J}^*\colon{H}^2_\chi\longrightarrow\Gamma}$ by the following change of orthonormal bases
\[\mathcal{J}\colon\Gamma\ni\mathfrak{e}^{\odot \lambda}_\imath\left\|\mathfrak{e}^{\odot\lambda}_\imath\right\|^{-1}
\longmapsto\varepsilon^\lambda_\imath\left\|\varepsilon^\lambda_\imath\right\|^{-1}_{L^2_\chi}\in{H}^2_\chi,
\qquad \lambda\in\mathbb{Y},\quad \imath\in\mathbb{N}^{\ell(\lambda)}_\ast.
\]
Clearly, $\mathcal{J}^*\colon\varepsilon^\lambda_\imath\left\|\varepsilon^\lambda_\imath\right\|^{-1}_{L^2_\chi}
\longmapsto\mathfrak{e}^{\odot \lambda}_\imath\left\|\mathfrak{e}^{\odot\lambda}_\imath\right\|^{-1}$, because
$\left\langle\mathcal{J}\mathfrak{e}^{\odot\lambda}_\imath\mathrel{\big|} f\right\rangle_{\!{L^2_\chi}}=
\left\langle\mathfrak{e}^{\odot\lambda}_\imath\mathrel{\big|} \mathcal{J}^*f\right\rangle$ for any ${f\in H^2_\chi}$.
Using Theorem~\ref{irrep1}, for any element ${\psi\in\Gamma}$ with the Fourier coefficients $\hat\psi_{(\lambda,\imath)}=\langle\mathfrak{e}^{\odot\lambda}_\imath\mid\psi\rangle
\|\mathfrak{e}^{\odot\lambda}_\imath\|^{-2}$, we obtain
\[
\mathcal{J}\psi=\sum_{(\lambda,\imath)\in\mathbb{Y}\times\mathbb{N}^{\ell(\lambda)}_\ast}
\hat\psi_{(\lambda,\imath)}
\frac{\|\mathfrak{e}^{\odot\lambda}_\imath\|^2}{\|\varepsilon^\lambda_\imath\|^2_{L^2_\chi}}
\varepsilon^\lambda_\imath\quad\text{where}\quad
\frac{\|\mathfrak{e}^{\odot\lambda}_\imath\|^2}{\|\varepsilon^\lambda_\imath\|^2_{L^2_\chi}}=
\frac{(\ell(\lambda)-1+|\lambda|)!}{(\ell(\lambda)-1)! |\lambda|!}.
\]
In particular, $\mathcal{J}x=\sum\hat{x}_k\varepsilon_k$ for any elements ${x\in\mathsf{E}}$ with
the Fourier coefficients $\hat{x}_k={\langle x\mid\mathfrak{e}_k\rangle}$. Moreover,
$\|\mathcal{J}x\|_{L^2_\chi}^2=\sum\|\hat{x}_k\|^2=\|x\|^2$.

In what follows, we assign to each $x\in\mathsf{E}$ the $L^2_\chi$-valued function
\[
x_\mathcal{J}\colon\mathfrak{U}\ni u\longmapsto(\mathcal{J}x)(u).
\]

\begin{lem}\label{infty}
The function $\mathcal{J}(1-x)^{-\otimes1}=(1-x_\mathcal{J})^{-1}$ in the variable $u\in\mathfrak{U}$
takes values in  $L_\chi^2$ for all ${x\in\mathsf{B}}$
\end{lem}

\begin{proof}
Applying $\mathcal{J}$ to the decompositions \eqref{Fx} and \eqref{Tayl}, we obtain
\begin{equation}\label{7}
\begin{split}
\mathcal{J}(1-x)^{-\otimes1}&=\sum_{(\lambda,\imath)\in\mathbb{Y}\times\mathbb{N}^{\ell(\lambda)}_\ast}
\frac{\hat{x}^\lambda_\imath\varepsilon^\lambda_\imath}{\|\mathfrak{e}^{\odot\lambda}_\imath\|^2}\\
&=\sum_{n\in\mathbb{Z}_+}\Big(\sum_{k\in\mathbb{N}}\hat{x}_k\varepsilon_k\Big)^n=(1-x_\mathcal{J})^{-1}
\end{split}
\end{equation}
where  the following orthogonal series with a fixed $n\in\mathbb{N}$,
\begin{equation}\label{8}
x_\mathcal{J}^n=\Big(\sum_{k\in\mathbb{N}}\hat{x}_k\varepsilon_k\Big)^n=
\sum_{(\lambda,\imath)\in\mathbb{Y}_n\times\mathbb{N}^{\ell(\lambda)}_\ast}
\frac{\hat{x}^\lambda_\imath\varepsilon^\lambda_\imath}{\|\mathfrak{e}^{\odot\lambda}_\imath\|^2},
\end{equation}
 is convergent  in $L^2_\chi$. Moreover, taking into account  the orthogonality, we get
\begin{align*}
\left\|(1-x_\mathcal{J})^{-1}\right\|_{L^2_\chi}^2&=\sum_{n\in\mathbb{Z}_+}
\sum_{(\lambda,\imath)\in\mathbb{Y}_n\times\mathbb{N}^{\ell(\lambda)}_\ast}
\frac{|\hat{x}^\lambda_\imath|^2}{\|\mathfrak{e}^{\odot\lambda}_\imath\|^2}\\
&=\sum_{n\in\mathbb{Z}_+}\Big(\sum_{k\in\mathbb{N}}|\hat{x}_k|^2\Big)^n=\left(1-\|x\|^2\right)^{-1}.
\end{align*}
Hence, the function $(1-x_\mathcal{J})^{-1}$ with $x\in\mathsf{B}$ takes values in  $L_\chi^2$.
\end{proof}

Let $f={\sum_n f_n\in{H}^2_\chi}$ with $f_n\in H_\chi^{2,n}$. Then ${\mathcal{J}^* f\in\Gamma}$ and
${\mathcal{J}^* f_n\in\mathsf{E}^{\odot n}}$.
Briefly denote $\tilde{f}:=(\mathcal{J}^* f)^*\in H_n^2$ and $\tilde{f}_n:=(\mathcal{J}^* f_n)^*\in H^2$. Thus,
\begin{align*}
\tilde{f}(x)&=\left\langle(1-x)^{-\otimes1}\mid\mathcal{J}^* f\right\rangle,\qquad x\in\mathsf{B},\\
\tilde{f}_n(x)&=\left\langle{x}^{\otimes n}\mid\mathcal{J}^* f_n\right\rangle,\qquad x\in\mathsf{E}.
\end{align*}

\begin{thm}\label{hard3}
Each Hilbert-Schmidt analytic function $\tilde{f}\in H^2$ has the integral representation
\begin{equation}\label{laplaceA}
\tilde{f}(x)=\int_\mathfrak{U}\frac{f\,d\chi}{1-x_\mathcal{J}},\qquad x\in\mathsf{B}
\end{equation}
and its Taylor coefficients at zero have the form
\begin{equation}\label{TaylorL}
\frac{d^n_0\tilde{f}(x)}{n!}=\int_\mathfrak{U}x_\mathcal{J}^nf_n\,d\chi,\qquad x\in\mathsf{E}.
\end{equation}
The mapping $f\longmapsto \tilde{f}$ produces the linear isometry ${H}^2_\chi\simeq{H}^2$.
\end{thm}

\begin{proof}
Consider the Fourier decomposition of $f$ with respect to the basis $\varepsilon^\mathbb{Y}$ and its $\mathcal{J}^*$-image, respectively
\[
f=\sum_{(\lambda,\imath)\in\mathbb{Y}\times\mathbb{N}^{\ell(\lambda)}_\ast}
\hat{f}_{(\lambda,\imath)}{\varepsilon^\lambda_\imath},\qquad
\mathcal{J}^*f=\sum_{(\lambda,\imath)\in\mathbb{Y}\times\mathbb{N}^{\ell(\lambda)}_\ast}
\bar{\hat{f}}_{(\lambda,\imath)}
\frac{\|\varepsilon^\lambda_\imath\|^2_{L^2_\chi}}{\|\mathfrak{e}^{\odot\lambda}_\imath\|^2}
\mathfrak{e}^{\odot\lambda}_\imath
\]
where $\hat{f}_{(\lambda,\imath)}=\|\varepsilon^{\lambda}_\imath\|_{L^2_\chi}^{-2}
\int_\mathfrak{U}{f}\,\bar\varepsilon^\lambda_\imath\,d\chi$.
Substituting $\hat{f}_{(\lambda,\imath)}$ to $\tilde{f}=(\mathcal{J}^* f)^*$ and using the orthogonal property and the relations
 \eqref{Tayl} and  \eqref{7}, we obtain
\begin{align*}
\tilde{f}(x)&
=\sum_{(\lambda,\imath)\in\mathbb{Y}\times\mathbb{N}^{\ell(\lambda)}_\ast}
\frac{\hat{f}_{(\lambda,\imath)}
\hat{x}^\lambda_\imath\left\langle \mathfrak{e}^{\odot\lambda}_\imath\mid\mathfrak{e}^{\odot\lambda}_\imath\right\rangle
\|\varepsilon^\lambda_\imath\|^2_{L^2_\chi}}{\|\mathfrak{e}^{\odot\lambda}_\imath\|^4}\\
&=\int_\mathfrak{U}\sum_{(\lambda,\imath)\in\mathbb{Y}\times\mathbb{N}^{\ell(\lambda)}_\ast}
\frac{\hat{x}^\lambda_\imath\varepsilon^{\lambda}_\imath}{\|\mathfrak{e}^{\odot\lambda}_\imath\|^2}f\,d\chi
=\int_\mathfrak{U}\frac{f\,d\chi}{1-x_\mathcal{J}}.
\end{align*}
Hence, \eqref{laplaceA} holds. Using \eqref{8}, we similarly obtain
\begin{equation}\label{n}
\tilde{f}_n(x)=\left\langle{x}^{\otimes n}\mathrel{\big|} \mathcal{J}^*f_n\right\rangle
=\int_\mathfrak{U}x_\mathcal{J}^nf_n\,d\chi.
\end{equation}
Taking into account \eqref{n} and the orthogonal decomposition \eqref{ort}, we get
\begin{equation}\label{r}
\tilde{f}\left(\alpha{x}\right)=
\left\langle(1-\alpha{x})^{-\otimes1}\mathrel{\big|}\mathcal{J}^*f\right\rangle=
\sum\alpha^n\int_\mathfrak{U}x_\mathcal{J}^nf_nd\chi,\quad
{|\alpha|\le1}.
\end{equation}
Note that $\tilde{f}\left(\alpha{x}\right)$ is analytic in $\alpha$ for all ${x\in\mathsf{B}}$.
Differentiating $\tilde{f}\left(\alpha{x}\right)$ at $\alpha=0$ and using the $n$-homogeneity of derivatives, we obtain
\[
\frac{d^n}{d\alpha^n}
\sum\alpha^n\int_\mathfrak{U}x_\mathcal{J}^n{f}_n\,d\chi
\mathrel{\Big|}_{\alpha=0}=n!\int_\mathfrak{U}x_\mathcal{J}^n{f}_n\,d\chi.
\]
Hence, the functions \eqref{TaylorL} coincide with the Taylor coefficients at zero of  $\tilde{f}$.

Finally, since the image of $\varepsilon^\mathbb{Y}$ under $\mathcal{J}^*$
coincides with $\mathfrak{e}^{\odot\mathbb{Y}}$, the mapping $H^2_\chi\ni f\longmapsto \tilde{f}\in H^2$
is an isometry.
\end{proof}

\section{Radial boundary values}

Using \eqref{laplaceA},  for each $f={\sum_n f_n\in{H}^2_\chi}$ with $f_n\in H_\chi^{2,n}$  we can rewrite \eqref{r} as
\[
\tilde{f}(rx)=\left\langle(1-r{x})^{-\otimes1}\mathrel{\big|}\mathcal{J}^*f\right\rangle
=\int_\mathfrak{U}\frac{f\,d\chi}{1-rx_\mathcal{J}},\qquad x\in\mathsf{K},\quad {r\in[0,1)}
\]
where $\mathsf{K}=\left\{x\in\mathsf{E}\colon\|x\|\le1\right\}$.
\begin{thm}\label{car:hardy2}
The integral transform $\mathcal{C}_r\colon{f}\longmapsto \mathcal{C}_r[f]$, defined as
\begin{equation}\label{CauchyB}
\mathcal{C}_r[f](x):=\int_\mathfrak{U}\frac{f\,d\chi}{1-rx_\mathcal{J}},\qquad x\in\mathsf{K},\quad{r\in[0,1)},
\end{equation}
belongs to the space of bounded linear operators $\mathscr{L}(H^2_\chi,H^2)$.
The radial  boundary values of ${\mathcal{C}_r[f]\in H^2}$
are equal to $\tilde{f}\in   H^2$ in the following sense:
\begin{equation}\label{boudval}
\lim_{r\nearrow1}\big\|\mathcal{C}_r[f]-\tilde{f}\big\|_{H^2}=0.
\end{equation}
Moreover, the following equality  holds,
\begin{equation}\label{h2norm}
\|\tilde{f}\|_{H^2}^2=\sup_{r\in[0,1)}\left\|\mathcal{C}_r[f]\right\|^2_{H^2}.
\end{equation}
\end{thm}

\begin{proof} Theorem~\ref{hard3} and \eqref{CauchyB} imply the equality
$\mathcal{C}_r[f]=\sum r^n\tilde{f}_n$ for any ${r\in[0,1)}$.
By   \eqref{iso}, we have
 $\tilde{f}_k\perp\tilde{f}_n$  as $n\noteq k$ in $H^2$. It follows that
\[
\left\|\mathcal{C}_r[f]\right\|^2_{H^2}=\left\|\sum r^n\tilde{f}_n\right\|^2_{H^2}
=\sum r^{2n}\|\tilde{f}_n\|^2_{H^2}= \sum r^{2n}\|f_n\|^2_{L_\chi^2},
\]
since $\mathcal{J}^*$ acts isometrically from $H^{2,n}_\chi$ onto the space $\mathsf{E}^{\odot n}$ which is
antilinear isometric to $H^2_n$ by definition. Similarly, we obtain that
\[
\big\|\mathcal{C}_r[f]-\tilde{f}\big\|^2_{H^2}=
\sum\left(r^{2n}-1\right)\|f_n\|^2_{L_\chi^2}\longrightarrow0,\qquad
r\to1.
\]
Moreover, the Cauchy-Schwarz inequality implies that
\[
\left\|\mathcal{C}_r[f]\right\|^2_{H^2}
\le\frac{1}{(1-r^2)^{1/2}}\Big(\sum\left\|f_n\right\|^2_{L_\chi^2}\Big)^{1/2}
=\frac{\|f\|_{L^2_\chi}}{(1-r^2)^{1/2}}
\]
for all $ f\in H^2_\chi$.
Hence, the operator $\mathcal{C}_r$ belongs to $\mathscr{L}(H^2_\chi,H^2)$ for all  $r\in [0,1)$.

Finally, the equalities
\[
\sup_{r\in[0,1)}\left\|\mathcal{C}_r[f]\right\|^2_{H^2}=\sup_{r\in [0,1)}\sum r^{2n}\|\tilde{f}_n\|^2_{H^2}
=\sum\|\tilde{f}_n\|^2_{H^2}=\|\tilde{f}\|^2_{H^2}
\]
give the required formula \eqref{h2norm}.
\end{proof}

%\subsection*{Acknowledgment} I am grateful to the referees for their comments and valuable suggestions which greatly improved this article.

% ------------------------------------------------------------------------

\begin{thebibliography}{1}
\bibitem{BerezanskiKondratiev95}
{Yu. M. Berezanski  \and Yu.G.  Kondratiev},
\textit{Spectral methods in infinite-dimensional analysis}.
Springer, 1995.

\bibitem{BorodinOlshanski05}
A. Borodin \and G. Olshanski,  \textit{Harmonic analysis on the
  infinite-dimensional unitary group and determinantal point
  processes}. Ann. Math. \textbf{161} (2005), 1319--1422.

\bibitem{Borodin11}
A. Borodin, \textit{Determinantal point processes} in \textit{Oxford Handbook of
Random Matrix Theory} (G. Akemann, J. Baik, and P. Di Francesco,
eds.) Oxford Univ. Press, 2011.

\bibitem{Carado09}
D. Carado, V. Dimand \and  S. Muro,
\textit{Coherent sequences of polynomial ideals on Banach spaces}.
Math. Nachr. \textbf{282}(8) (2009), 1111--1133.

\bibitem{ColeGamelin86}
B. Cole \and  T.W.  Gamelin,
\textit{Representing measures and {Hardy} spaces for the infinite polydisk  algebra}.
Proc. London Math. Soc. \textbf{53} (1986), 112--142.

\bibitem{Dwyer71}
T.A.W. Dwyer III,
\textit{Partial differential equations in Fischer-Fock spaces for the Hilbert-Schmidt holomorphy type}.
{Bull. Amer. Math. Soc.} \textbf{77}(5) (1971), 725--739.

\bibitem{G}
T.W. Gamelin, \textit{Analytic functions on {Banach} spaces},  in
\textit{Complex Function  Theory}  (Gauthier and Sabidussi eds.)
Kluwer, 1994,  187--223.

\bibitem{Floret97}
K. Floret, \textit{Natural norms on symmetric tensor products of normed spaces}.
Note di Matematica \textbf{17},  (1997), 153--188.

\bibitem{Herv}
M.~Herv\'e, \emph{Analyticity in {Infiite} {Dimensional}
{Spaces}}, de Gruyter  Stud. in Math., vol.‾10, Walter de Gruyter, Berlin, New York, 1989.

\bibitem{HewittRoss70}
{E. Hewitt \and K. A.  Ross},
\textit{Abstract Harmonic Analysis}, Vol.2, Springer, 1994.

\bibitem{lopushansky2013}
O. Lopushansky,
\textit{Hardy type space associated with an infinite-dimensional unitary matrix group}.
Abst. Appl. An. ID 810735 (2013),  1--7.


\bibitem{Neretin2002}
Yu. A. Neretin,  \textit{Hua type integrals over unitary groups and over projective limits
of unitary groups}. Duke Math. J. \textbf{114}(2) (2002), 239--266.


\bibitem{Olshanski2003}
G. Olshanski, \textit{The problem of harmonic analysis on the
infinite-dimensional unitary group}.
J. Funct. Analysis. \textbf{205} (2003), 464--524.

\bibitem{OrtedNeeb98}
{K.H. Neeb \and B. {\O}rted},
\textit{Hardy spaces in an infinite dimensional setting},
In: {H.D.Doebner (ed.)},
\textit{Lie Theory and Its Applications in Physics,  H.D.Doebner (ed.)},
{3 -- 27}, {Word Sci. Publ.}, 1998.

\bibitem{Petersson2001}
P. Petersson, \textit{Hypercyclic convolution operators on entire functions of
Hilbert-Schmidt holomorphy type}.
Ann. Math. Blaise Pascal \textbf{8}(2) (2001), 107--114.

\bibitem{Pickrell}
D. Pickrell, \textit{Measures on infinite-dimensional Grassmann
 manifolds}. J. Funct. Analysis. \textbf{70} (1987), 323--356.

\bibitem{PinascoZalduendo05}
D. Pinasco \and  I. Zalduendo,
\textit{Integral representations of holomorphic functions on {Banach} spaces},
J. Math. Anal. Appl. \textbf{308} (2005), 159--174.

\bibitem{RudinFT80}
W. Rudin,  \textit{Function {theory} in the {unit} {ball} of $\mathbb{C}^n$}.
Springer, 2008.

\bibitem{Sait}
S.Saitoh, \textit{Integral Transforms, Reproducing Kernels and Their
Applications}. Pitman Research Notes in Math. Ser. Vol.~369,
Longman, 1997.
\end{thebibliography}
\end{document}